\documentclass[letterpaper,12pt,leqno,oneside]{article}
\usepackage{color}
\usepackage{bm}
\usepackage{exscale}
\usepackage{amsmath}
\usepackage{amsfonts}
\usepackage{wasysym}
\usepackage{stmaryrd}
\usepackage{amscd}
\usepackage{graphicx}
\usepackage{pagecolor,lipsum}
\usepackage{xcolor}
\usepackage{textcomp}
\usepackage{amsxtra}
\usepackage{amssymb}
\usepackage{theorem}
\usepackage[final]{epsfig}
\usepackage{eqnarray}
\usepackage{hyperref}
\newtheorem{proposition}{Proposition}[subsection]
\newtheorem{definition}[proposition]{Definition}

\newtheorem{lemma}[proposition]{Lemma}

{\theorembodyfont{\rmfamily}}

{\theorembodyfont{\rmfamily}}

\newfont{\abc}{cmtt10 scaled 1200}

\def\R{\mathbb{R}}

\def\Z{\mathbb{Z}}

\def\P{\mathbb{P}}

\def\P{\mathbb{P}}

\def\ve{\varepsilon}

\def\ve{\varepsilon}

\def\ra{\rightarrow}

\def\p{\partial}

\def\qed{\hfill $\Box$ \\}
\def\mm{\mbox}
\def\v{= \emptyset}
\def\n{\neq \emptyset}

\def\bp{\langle A \rangle}

\def\si{$\mathcal{S}$}

\begin{document}

\vspace*{-0.3cm}

\begin{center}\Large{\bf{Minimal Smoothings of Area Minimizing Cones}}\\
\bigskip
\large{\bf{Joachim Lohkamp}\\
\medskip}

\end{center}

\noindent Mathematisches Institut, Universit\"at M\"unster, Einsteinstrasse 62, Germany\\
 {\small{\emph{e-mail: j.lohkamp@uni-muenster.de}}}

{\small {\center \tableofcontents}

\setcounter{section}{1}
\renewcommand{\thesubsection}{\thesection}
\subsection{Introduction} \label{int}
Area minimizing hypersurfaces may contain a complicated singularity set $\Sigma\subset H$ and $H \setminus \Sigma$ degenerates towards $\Sigma$ in a rather delicate way. It is well-known that the codimension
of $\Sigma$ is $\ge 7$. This estimate and the entire regularity theory of such hypersurfaces follows from inductive use and analysis of again area minimizing tangent cones we get around each singular point of $H$.
This is the basic example of the quite common strategy of cone reduction arguments in this area which are the major reason for the interest in Euclidean area minimizing cones.\\

In this paper we show that every area minimizing cone $C^n \subset \R^{n+1}$ can be approximated by entirely \textbf{smooth}  area minimizing hypersurfaces. This was previously known only when $C^n$ admits at most an isolated singularity in $0$ by Hardt and Simon \cite{HS}, Th.2.1. Their work was further refined by Macintosh \cite{Mc} and Mazzeo and Smale  \cite{MS}.  It is also worthy to mention that, for the Simons cone, the Hardt-Simon result is implicity contained in the earlier work of Bombieri, DeGiorgi and Giusti \cite{BDG} when they settled the Bernstein problem.\\

To formally state our main result we recall that any area minimizing cone $C^n \subset \R^{n+1}$ separates its complement into two  path components $E_C^+ \cup E_C^-=\R^{n+1} \setminus C$, cf.\cite{BG}.
Also, we recall the notion of an oriented minimal boundary. This is a Euclidean area minimizing hypersurface we can write as the boundary of an open subset of $\R^{n+1}$. We refer to \cite{L1},Appendix A for a short survey covering those details needed in our context.\\

We note that either $C^n \subset \R^{n+1}$  is singular at least in $0$ and this may happen in dimension $n+1 \ge 8$ or  $C^n \subset \R^{n+1}$ is a smooth hyperplane.\\

\textbf{Theorem} \textbf{(Minimal Cone Smoothing)} {\itshape \,  Let $C^n \subset \R^{n+1}$ be an area minimizing cone and $E_C$ be one of the two open components $E_C^\pm$ of $\R^n \setminus C$. Then we have the following smoothing result:
\begin{itemize}
  \item \emph{\textbf{Regularity}} \, There is a \textbf{smooth} oriented minimal boundary $H_C \subset E_C$ and  $C$  is its unique tangent cone at infinity.
  \item \emph{\textbf{Controlled Foliations}} \,The hypersurfaces $s \cdot H_C$, $s>0$  foliate $E_C$ and they exclusively induce Jacobi fields* on $C$ having \textbf{minimal growth} towards the singularities of $C$.
  \item \emph{\textbf{Uniqueness}} \, For any  oriented minimal boundary $T \subset E_C$, there is some constant $c_T >0$ so that $T = c_T \cdot H_C$.
\end{itemize}}

*Throughout the paper a \textbf{Jacobi field}  means a function $f$ solving $J_H \, w =0$, for the Jacobi field operator $J_H=-\Delta_H - |A|^2$. Here $|A|$ is the norm of the second fundamental form $A$ of $H \subset \R^{n+1}$. Of course, $f$ is just the size of the normal component of a Jacobi (vector) field.\\

\textbf{Outline of the Proof} \, Our argument  differs substantially from that in \cite{HS}. It
is based on Martin theory of the Jacobi field operator on these cones we approach by means of \textbf{hyperbolic unfoldings} of these cones \cite{L1} - \cite{L3}. This makes the cone smoothing more transparent since we drop the analysis of the nonlinear asymptotic presentation of area minimizers towards infinity which was an essential part of the argument in \cite{HS} and which limited its applicability to isolated singularities. Turning to some of the details, we first notice that the existence assertion of a possibly singular hypersurface $H_C \subset E_C$ is rather simple. The non-trivial part is its smoothness. The proof is by induction starting from cones with isolated singularities which show up in all dimensions $\ge 8$. For this short overview we focus on cones singular in exactly $n-7$ dimensions. \\

 In dimension $8$ we consider some area minimizing cone $C^7 \subset \R^8$ and show that $H_{C^7} \subset E_{C^7}$ induces an, in fact uniquely determined, Jacobi field on $C$ with minimal growth towards $0$.   Since we have a detailed picture of such Jacobi fields we can give a simple argument to show the smoothness and also the uniqueness of $H_{C^7} \subset E_{C^7}$ (up to scalings). To prepare the induction step, we also show that the smooth hypersurface $H_{C^7}\times \R \subset E_{{C^7} \times \R}$ is again the unique such area minimizer in $E_{{C^7} \times \R}$ and that it induces the unique Jacobi field on ${C^7} \times \R$ with minimal growth towards $\{0\}\times \R$.\\

Then we turn to the next higher dimension. We consider some minimizing cone $C^8 \subset \R^9$ and some associated area minimizing  hypersurface $H_{C^8} \subset E_{C^8}$. The main point is to show that $H_{C^8}$ induces a Jacobi field on $C^8$ with minimal growth towards
the singular set $\Sigma_C$ of $C^8$. Here we use that we have already derived this result for all tangent cones of  $C^8$ outside $0$ and infer that the Jacobi field on $C^8$ must have minimal growth towards $\Sigma_C \setminus \{0\}$. Then we further show that this also holds in $\{0\}$. From this we proceed as before. We show the smoothness and uniqueness of $H_{C^8} \subset E_{C^8}$ and continue inductively.\\

The potential theory we import from \cite{L1} - \cite{L3} controls all the correlations between the Jacobi fields which have minimal growth towards the singular set on the given cone, on its tangent cones and also the approximating smooth hypersurfaces.\\

\textbf{Applications} \, An immediate consequence of the cone smoothing theorem is that, since any area minimizing hypersurface $H^n$, in a Riemannian manifold $(M^{n+1},g)$ can locally be approximated from a tangent cone, it can equally be approximated by portions of smooth area minimizers. At first, these area minimizers will belong to  $\R^{n+1}$, but they can be perturbed to (small pieces of) area minimizers in $(M^{n+1},g)$. This supports the classical \textbf{smoothing conjecture} that, after small perturbations of $g$ in $C^k$-topology, for any $k \ge0$, one may assume that $H^n$ is smooth. In fact, Smale \cite{Sm} has made this strategy work in the case of isolated singularities. For higher dimensional singularities it is yet unclear how to use the, then interfering, local smoothings to prove the smoothing conjecture.\\

 However, the cone smoothing already implies a weaker smoothing result also conjectured since the late 70ties: the  \textbf{splitting theorem in scalar curvature geometry}. It can interpreted as a counterpart of the Cheeger-Gromoll splitting theorem in Ricci curvature geometry:\\

 {\itshape \,  Let $(M^{n+1},g)$, $n \ge 3$, be a smooth, compact Riemannian manifold with positive scalar curvature  and $\alpha \in H_n(M;\Z)$. Then, there is a \textbf{smooth} compact hypersurface $H^n\subset M^{n+1}$ that represents $\alpha$ and admits a \textbf{smooth}  positive scalar curvature metric.}\\

From Smale's approach in  \cite{Sm}, establishing the smoothing conjecture for isolated singularities, one readily gets this splitting result in dimension $8$ where singularities are always isolated.
The idea for the general case is to simply exchange the order of operations in this argument. Again, we start from a (usually singular) area minimizing hypersurface $N^n \in \alpha$. But, then we first conformally deform $N^n$, its tangent cones and the approximating smooth hypersurfaces of these cones to matching positive scalar curvature geometries. Then, in a second step, positivity of scalar curvature is an open differential condition we can use to build a \textbf{smooth} positive scalar curvature geometry $H^n$ replacing singular patches of $N^n$ for smooth one.  This entails further details from the potential theory of the conformal Laplacian in \cite{L3}, again approachable via hyperbolic unfoldings and will be presented elsewhere.

\setcounter{section}{2}
\renewcommand{\thesubsection}{\thesection}
\subsection{Asymptotic Analysis of Jacobi Fields} \label{2}

The technical main ingredients for the proof of the cone smoothing theorem come from the potential theory of the Jacobi field operator. This potential theory is best described using so-called \si-structures and hyperbolic unfoldings, cf. \cite{L1} - \cite{L3} for details. In this chapter give a condensed introduction to this theory to draw some customized consequences needed in our later argument.

\subsubsection{Hyperbolic Unfoldings} \label{a1}

Let $M^{n+1}$ be a smooth  manifold, and $H^n\subset M^{n+1}$ be an area minimizing (or more generally almost minimizing) hypersurface with singularity set $\Sigma\subset H$. It is a compact set of Hausdorff-dimension $\le n-7$ and $H$ degenerates towards $\Sigma$. We shall consider the following classes of almost minimizers:
\begin{description}
  \item[${\cal{G}}^c_n$:] $H^n \subset M^{n+1}$ is a compact embedded \emph{almost} minimizer. We set ${\cal{G}}^c :=\bigcup_{n \ge 1} {\cal{G}}^c_n.$,
  \item[${\cal{H}}^{\R}_n$:] $H^n \subset\R^{n+1}$ is oriented minimal boundary in flat Euclidean space $(\R^{n+1},g_{eucl})$ with $0\in H$,
\item[${\cal{G}}_n$:] ${\cal{G}}_n := {\cal{G}}^c_n \cup {\cal{H}}^{\R}_n$ and ${\cal{G}} :=\bigcup_{n \ge 1} {\cal{G}}_n$. The point about  ${\cal{G}}_n$ is that it is closed under arbitrary blow-ups around singular points of elements $H \in {\cal{G}}_n$.
\end{description}
We notice that for any area minimizing cone $C^n \subset\R^{n+1}$, the hypersurface  $\p B_1(0) \cap C^n \subset S^n$ is \emph{not} area minimizing but it is an almost minimizer.\\

On these hypersurfaces we introduce natural distance and size concepts, the \emph{\si-structures}, which measure also the curvature of $H$. The basic concept is that of a \si-transform and we are using the more specific Hardy \si-transform we informally describe as follows.\\

\begin{definition} We call an assignment $\bp$ associating to any $H \in {\cal{G}}$ a  non-negative, measurable function $\bp_H:H \setminus \Sigma_H\to\R$ a  \textbf{Hardy \si-transform} provided it satisfies the following axioms:\\

$\bullet$\,  $\bp_H$  is naturally assigned to  $H$, that is, the assignment commutes
with the convergence of sequences of underlying area minimizers.\\

$\bullet$\, If $H$ is totally geodesic  $\bp_H \equiv 0$. Otherwise the level sets of $|A|$:
  \[
  \bp_H>0, \bp_H \ge |A_H| \mm{ and } \bp_H(x) \ra \infty, \mm{ for } x \ra p \in \Sigma_H.
  \]
Like $|A_H|$, $\bp_H$ anticommutes with scalings, i.e., $\bp_{\lambda \cdot H} \equiv \lambda^{-1} \cdot  \bp_{H}$ for any $\lambda >0$.\\

$\bullet$\, For $H \in{\cal{H}}^{\R}_n$  there exists a constant $k_n > 0$ depending only on the dimension such that  $f \in C_0^\infty(H \setminus \Sigma)$ we have the Hardy type inequality:  $\int_H|\nabla f|^2 + |A|^2 \cdot f^2 dV \ge k_n  \cdot \int_H \bp^2 \cdot f^2 dV$.\\

 $\bullet$\, If $H$ is not totally geodesic, and thus $\bp_H>0$, we define the $ \mm{ \textbf{\si-distance} } \delta_{\bp_H}:=1/\bp_H.$
  This function is
  $L_{\bp}$-Lipschitz regular for some constant $L_{\bp}=L(\bp,n)>0$, i.e.,
  \[
  |\delta_{\bp_H}(p)- \delta_{\bp_H}(q)|   \le L_{\bp} \cdot d_{g_H}(p,q) \mm{ for any } p,q \in  H \setminus \Sigma \mm{ and any } H \in {\cal{G}}_n.
  \]
\end{definition}
*The existence  of such a Hardy \si-transform and a broader explanation of these axioms is described in \cite{L2}.  For the remainder of this paper we choose one fixed Hardy \si-transform.

\begin{proposition} For any non-totally geodesic hypersurface $H \in {\cal{G}}$, and any pair $x$, $y\in H \setminus \Sigma$ we define  the \textbf{\si-metric}  $
d_{\bp}(x,y) := \inf   \{\int_\gamma  \bp \, \, \Big| \, \gamma   \subset  H \setminus \Sigma\mbox{ rectifiable curve joining }  x \mbox{ and } y  \}.$\\

$\bullet$\,  The Lipschitz regular metric space $(H \setminus \Sigma, d_{\bp}) = (H \setminus \Sigma, 1/\delta_{\bp}^2 \cdot g_H)$ and its quasi-isometric Whitney smoothing, i.e.\ the smooth Riemannian manifold $(H \setminus \Sigma, d_{\bp^*}) = (H \setminus \Sigma, 1/\delta_{\bp^*}^2 \cdot g_H)$, are \textbf{complete Gromov hyperbolic spaces} with \textbf{bounded geometry}. \\

$\bullet$\,  We refer to both these spaces as \textbf{hyperbolic unfoldings} of $(H \setminus \Sigma, g_H)$.\\

$d_{\bp}$ is natural. That is, the assignment of $d_{\bp_H}$ to $H$ commutes with the compact convergence of the regular portions of the underlying area minimizers.\\

$\bullet$\,  The identity map on $H \setminus \Sigma$ extends to homeomorphisms between the one-point compactification* $\widehat{H}$ and the Gromov compactifications of $(H \setminus \Sigma,d_{\bp})$ and $(H \setminus \Sigma,d_{\bp^*})$ so that $\widehat{H}\cong\overline{(H \setminus \Sigma,d_{\bp})}_G \cong \overline{(H \setminus \Sigma,d_{\bp^*})}_G$,
where $ \cong$ means homeomorphic. Thus, we find for the associated \textbf{Gromov boundaries} $\p_G$: $\widehat{\Sigma} \cong\p_G(H \setminus \Sigma,d_{\bp}) \cong \p_G(H \setminus \Sigma,d_{\bp^*})$.
\end{proposition}

\noindent\textbf{Proof} \,   \cite{L1}, Theorem 1.11 and Theorem 1.13 \qed

*The \emph{one-point compactification} of a hypersurface $H \in {\cal{H}}^{\R}_n$ is written $\widehat{H}$. For the singular set $\Sigma_H$ of some $H \in {\cal{H}}^{\R}_n$ we \emph{always} add $\infty_H$  to $\Sigma$ and define $\widehat\Sigma:=\Sigma \cup \infty_H$ (note that $\Sigma$ could already be compact). For  $H \in {\cal{G}}^{c}_n$ we set $\widehat H=H$ and $\widehat\Sigma=\Sigma$.

\subsubsection{Martin  Theory and Critical Operators} \label{a2}
The point about hyperbolic unfoldings is that due to work of Ancona, \cite{KL} for an exposition, the potential theory of many uniformly elliptic second order operators on complete Gromov hyperbolic spaces with bounded geometry is as simple as the theory for the Laplace operator on the unit disc. For instance, their Gromov boundary is homeomorphic to the Martin boundary for these operators. A pull-back from the hyperbolic unfolding to the initial almost minimizer gives us a very transparent potential theory for the transformed elliptic operators near their singularities.\\

To describe the elliptic problems on $(H \setminus \Sigma, g_H)$ that we can address this way we use special charts for $H \setminus \Sigma$, namely \emph{\si-adapted charts}. These are bi-Lipschitz charts $\psi_p:B_R(p)\to\R$ centered in $p\in H\setminus\Sigma$ where the radius $R$ of the ball depends on $\bp_H(p)$, cf.\ \cite[B]{L1}.

\begin{definition} Let $H \in \cal{G}$. We call a symmetric second order elliptic operator $L$ on $H
\setminus \Sigma$  \textbf{shifted \si-adapted} supposed the following two conditions hold:\\

\textbf{$\mathbf{\bp}$-Adaptedness}  \,
$L$ satisfies \si-weighted uniformity conditions with respect to the charts $\psi_p$. Namely, we can write $-L(u) = \sum_{i,j} a_{ij} \cdot \frac{\p^2 u}{\p x_i \p x_j} + \sum_i b_i \cdot \frac{\p u}{\p x_i} + c \cdot u,$
\emph{for some locally} $\beta$-H\"{o}lder continuous coefficients $a_{ij}$, $\beta \in (0,1]$, measurable functions $b_i$ and $c$, and there exists a $k_L=k \ge 1$ such that for any $p \in H\setminus \Sigma$ and $\xi\in\R^n$:
\begin{itemize}
  \item $k^{-1} \cdot\sum_i \xi_i^2 \le \sum_{i,j} a_{ij}(p)\cdot \xi_i \xi_j \le k \cdot \sum_i \xi_i^2$,
  \item $\delta^{\beta}_{\bp}(p)\cdot |a_{ij}|_{C^\beta(B_{\theta(p)}(p))} \le k$,  $\delta_{\bp}(p)\cdot |b_i|_{L^\infty}\le k$ and $\delta^2_{\bp}(p) \cdot |c|_{L^\infty}\le k$.
\end{itemize}
\textbf{$\mathbf{\bp}$-Finiteness} \,
There is some $\tau = \tau(L,\bp,H)>-\infty$ such that for any smooth function $f$ which is compactly supported in $H\setminus \Sigma$, we have $\int_H  f  \cdot  L f  \,  dV \, \ge \, \tau \cdot \int_H \bp^2\cdot f^2 dV.$
The largest such $\tau$   is the \textbf{principal eigenvalue} $\lambda^{\bp}_{L,H}$ of $\delta_{\bp}^2 \cdot L$.  The operator $L$ is called \textbf{\si-adapted} if $\lambda^{\bp}_{L,H}>0$.
\end{definition}

In this paper, our example for such operators comes from the variation of minimal surfaces:

\begin{lemma}
For any  $H \in \cal{G}$, the \textbf{Jacobi field operator} $J_H=-\Delta_H - |A|^2-Ric_M(\nu,\nu)$ is shifted \si-adapted. Moreover, if $H \in \cal{H}$, then $J_H$ has principal eigenvalue $\ge 0$.
\end{lemma}
\noindent\textbf{Proof} \,   \cite{L3},Theorem 2. \qed

Our goal is to understand positive solutions of  $J_H \, w = 0$ on an area minimizing hypersurface $H$. One-sided variations of minimal hypersurfaces lead to Jacobi fields one may express in terms of such solutions pointwise multiplying the normal vector field of the hypersurface.   Minimal growth properties of such fields towards points in $\widehat{\Sigma}$ play a central role in the asymptotic analysis of $J_H$. On $H \setminus \Sigma$ a solution $u >0$ of the equation $J_H  \, f=0$ usually diverges to infinity when we approach $\Sigma$. The notion of \emph{$J_H$-vanishing} is a  \emph{minimal growth} condition generalizing that of  classically vanishing boundary data for the Laplacian on a disc.

\begin{definition}  A solution $u >0$ of $L \, f=0$ is said to be \textbf{$L$-vanishing} in $p \in \widehat{\Sigma}$, if there is a supersolution $w >0$, so that $u/w(x) \ra 0$, for $x \ra p$, $x \in H \setminus \Sigma$. We call $u >0$ \textbf{minimal} if for any other solution $v >0$ with $v \le u$, we have $v \equiv c \cdot u$ for some $c >0$. The space of minimal solutions (normalized to $1$ in some basepoint) is the \textbf{(minimal) Martin boundary}  $\p^0_M (H \setminus \Sigma,L)$.
\end{definition}

It can be shown that $L$-vanishing is equivalent to minimal growth and we occasionally prefer to call  $L$-vanishing solutions the solutions of minimal growth.  Now the central result says that, in our case, the usually very hard to determine Martin boundary is just the singular set.

\begin{proposition} \label{mc}
Let $H \in {\cal{G}}$ be not totally geodesic and $L$ an \si-adapted operator on $H \setminus \Sigma$. Then
\begin{itemize}
  \item the identity map on $H \setminus \Sigma$ extends to a homeomorphism  between $\widehat{H}$ and the \emph{\textbf{Martin compactification}} $\overline{H \setminus \Sigma}_M$.
  \item all Martin boundary points are minimal:
  $\p^0_M (H \setminus \Sigma,L) \equiv \p_M(H \setminus \Sigma,L)$.
\end{itemize}
Thus, $\widehat{\Sigma}$ and the minimal Martin boundary $\p^0_M (H \setminus \Sigma,L)$ are homeomorphic. A function $u >0$ on $H \setminus \Sigma$ solves $L \, f= 0$ iff there is a unique finite Radon measure $\mu$ on $\widehat{\Sigma}$ such that
\begin{equation}\label{muu}
u(x) = u_{\mu}(x) =\int_{\widehat{\Sigma}} k(x;y) \, d \mu(y).
\end{equation}

\end{proposition}

\noindent\textbf{Proof} \,   \cite{L2},Theorem 3. \qed

In this integral formula, $k(x;y)$ is the \emph{Martin kernel} of $L$ on $H \setminus \Sigma$. It is, up to multiples, the unique positive solution of $L \, f = 0$ on $H \setminus \Sigma$ which  $L$-vanishes in all points of $\widehat{\Sigma}$ except for $y$. Moreover, the functions $k(\cdot;y)$, $y\in\p_M(H \setminus \Sigma,L)$, are just the minimal solutions of $L$.\\

 If $L$ is  shifted \si-adapted on $H\setminus\Sigma$, we get the following basic spectral theory of $\delta_{\bp}^2 \cdot L$.

\begin{proposition} \label{tri}
Let $H \in \cal{G}$ and $\Sigma_H \n$, and $L$ be a shifted \si-adapted operator on $H \setminus \Sigma$. We set $L_\lambda:= L - \lambda \cdot \bp^2 \cdot Id,$ for $\lambda \in\R.$
Then we have the following trichotomy:
\begin{itemize}
    \item \emph{\textbf{Subcritical}} when $\lambda < \lambda^{\bp}_{L,H}$: $L_\lambda$ is \si-adapted.
    \item \emph{\textbf{Critical}} when $\lambda = \lambda^{\bp}_{L,H}$:  There exists, up to multiples, a unique positive solution $\phi$ of $L _{\lambda^{\bp}_{L,H}} \ f = 0$, the so-called ground state of $L_{\lambda^{\bp}_{L,H}}$.
    \item \emph{\textbf{Supercritical}} when $\lambda > \lambda^{\bp}_{L,H}$: $L_\lambda \, u = 0$ has no positive solution.
\end{itemize}
In the critical case $L_\lambda$ does not admit a positive Green's function. In turn, we say $L$ is \begin{itemize}
  \item \textbf{subcritical}, if $L$ admits a positive Green's function.
  \item \textbf{critical} if it does \emph{not} admit a Green's function but a positive solution of $L \, f=0$.
  \item \textbf{supercritical} when the latter equation does \emph{not} admit any positive solutions.\qed
\end{itemize}
\end{proposition}

\subsubsection{Applications to $J_C$} \label{a3}

The previously described Martin theory is again a natural theory. To explain this naturality we consider Schr\"odinger operators, like the Jacobi field operator, which are naturally associated to $H \in {\cal{G}}$ in the sense that the operators assigned to a converging sequence $H_i \in {\cal{G}}$ converge to that assigned to the limit hypersurface. Then also essential pieces of the potential theory commute under such limit processes. We discuss some of these results customized to our needs in the case of the Jacobi field operator $L=J$

\begin{proposition} \label{ver}
Let $H \in {\cal{H}}$ and $p \in \Sigma_H$ and $C$ be any tangent cone in $p$. Then we have $\lambda^{\bp}_{J,C} \ge \lambda^{\bp}_{J,H} \ge 0$ and the following two conditions are equivalent:\\
\begin{enumerate}
  \item A solution $u >0$ of $J_H \, f =0$ is $J_H$-vanishing in a neighborhood $V$ of $p$
  \item For any tangent cone $C^*$ in any point $q \in V$ and any solution of $J_{C^*} \, f =0$ induced by $u$ is $J_{C^*}$-vanishing along the entire singular set $\sigma_{C^*} \subset C^*$
\end{enumerate}
\end{proposition}

The result equally at the infinity and for \emph{tangent cones at infinity}.  Martin theory  shows that the induced solutions on each of these cones are uniquely determined and their associated Radon measure is the Dirac measure in $\infty$.\\

\noindent\textbf{Proof} \,  In \cite{L2},Theorem 4 we proved the main case of the inclusion $(i) \Rightarrow (ii)$ where $\lambda^{\bp}_{J,H} > 0$ and therefore $\lambda^{\bp}_{J,C^*} > 0$ on any tangent cone $C^*$.\\

When $\lambda^{\bp}_{J,H} = 0$, then $J$ is critical and there is only one positive solution, its ground state. It can be characterized as the limit of Dirichlet eigenfunctions on nested regular domains $H \setminus \Sigma_H$. This  shows it is $J_H$-vanishing towards every point on $\widehat{\Sigma}$.\\

Thus, if $\lambda^{\bp}_{J,H} = 0$ and $\lambda^{\bp}_{J,C^*} >  0$, we can apply the argument of \cite{L2},Prop.3.12, replacing the Green's function for the ground state, to conclude that the induced solutions are $J_{C^*}$-vanishing along the entire singular set $\sigma_{C^*} \subset C^*$. Otherwise, when $\lambda^{\bp}_{J,H} = 0$ and $\lambda^{\bp}_{J,C^*} = 0$, we have the ground state on $C^*$ and infer directly that it is $J_{C^*}$-vanishing along all of $\widehat{\sigma_{C^*}}$.\\

For the converse $(ii) \Rightarrow (i)$ we can assume that $\lambda^{\bp}_{J,H} = 0$, since the result is trivial when $\lambda^{\bp}_{J,H} = 0$. Now we argue by contradiction. That is, we assume we had a solution $u>0$ so that $u(x) = \int_{\widehat{\Sigma}} k(x;y) \, d \mu(y)$, for some Radon measure with $\mu(V)>0$ and so that any solution of $J_{C^*} \, f =0$ induced by $u$ is $J_{C^*}$-vanishing along the entire singular set $\sigma_{C^*} \subset C^*$, for any tangent $C^*$, in any point of $V$. Then we can choose a compact subset $K \subset V$ so that $\mu(K)>0$. \\

We inductively define nested covers ${\cal{C}}_k$  of $K$ by equisized balls with radii $2^{-k}$, for $k \ra \infty$. Since the ambient space is the $\R^{n+1}$ we can ensure upper bounded intersection numbers for these balls independent of $k$. Now we iteratively select nested balls $B_{2^{-k}}(p_k)$ carrying maximal measure $\mu(K \cap B_{2^{-k}}(p_k)$ under all such restrictions of $\mu$ onto the balls in ${\cal{C}}_k$.  The sequence then shrinks to a point $z \in K \subset V$. Now we turn to the hyperbolic unfolding  consider a neighborhood  basis defined from so-called $\Phi$-chains $\mathbf{N}^\delta_i(z) \subset H$, $i \in \Z^{\ge 0}$, of points $z \in \widehat{\Sigma}$ with  $\mathbf{N}^\delta_{i+1}(z)  \subset \mathbf{N}^\delta_i(z)$ and $\bigcap_i \mathbf{N}^\delta_i(z)=\{z\}$. We recall, that the subsets $\mathcal{N}^\delta_i(z) :=\mathbf{N}^\delta_i(z) \cap H \setminus \Sigma$ are nothing but halfspaces in the hyperbolic unfolding bounded by hypersurface perpendicular to a hyperbolic ray representing $z$
as a point in the Gromov boundary.\\

As in the argument for \cite{L2},Prop 3.15 Step 1 the $\Phi$-chains can be chosen to converge to $\Phi$-chains of the tangent cone $C^*$ in $z$. Now the argument of the inclusion $(i) \Rightarrow (ii)$ shows that the minimal solutions induce
minimal solutions in the limit and the convergence of the $\Phi$-chains ensures that the weight of the associated Dirac measure also converge. The choice of $z$ therefore implies that the Radon measure in the Martin integral of any induced solution on $C^*$ has strictly positive measure on $\Sigma_{C^*}$. This contradicts the $J_{C^*}$-vanishing of these solutions along $\Sigma_{C^*}$. \qed

The reason why we are particularly interested in $J_{C}$-vanishing solutions is that \ref{mc} also shows that they are uniquely determined on $C$, up to multiples.
This urges them to have very simple shape even in the presence of large singular portions on $C$ outside $0$. (By the way, while all this may appear fairly plausible, these are rather specific results for area minimizing cones. They become false for more general Euclidean cones. Ancona has given remarkable examples in \cite{A}.)

\begin{proposition} \label{ri} Let $C^n \subset \R^{n+1}$ a singular area minimizing cone. Then, we have $\lambda^{\bp}_{J,C} \ge 0$.\\

For $\lambda^{\bp}_{J,C} = 0$ the only positive solution of $J_C \, w=0$ is the ground state $\Psi$ and we have
\[
\Psi(\omega,r) = \psi(\omega) \cdot r^{- \frac{n-2}{2}}, \, (\omega,r) \in S_C \setminus \sigma \times \R^{>0}.
\]

For $\lambda^{\bp}_{J,C} > 0$ we have two distinguished points $\Psi_-$ at zero and $\Psi_+$ at infinity in the  Martin boundary of $J_C$, we have, in terms of polar coordinates $(\omega,r)$:
\[
\Psi_\pm(\omega,r) = \psi(\omega) \cdot r^{\alpha_\pm}, \, (\omega,r) \in S_C \setminus \sigma \times \R^{>0}, \mm{ with } \textstyle \alpha_\pm = - \frac{n-2}{2} \pm \sqrt{ \Big( \frac{n-2}{2} \Big)^2 + \mu}
\]
for some constant $0 >\mu(C)>-(n-2)^2/4$.
\end{proposition}

We observe that the restriction of the solutions $\Psi(\omega,r)$ and $\Psi_\pm(\omega,r)$, to any ray $\{\omega\} \times  \R^{>0} \subset C$, is an unbounded and strictly decreasing function in the variable $r>0$.\\

Finally we recall the following Fatou type theorem

\begin{proposition} \label{ft}
Let $H \in {\cal{G}}$ and $L$ be an \si-adapted operator on $H \setminus \Sigma$. Further, let $\mu$ and $\nu$ be two finite Radon measures on $\widehat{\Sigma}$ associated with solutions $u_{\mu}$ and $u_{\nu}$ of $L \, f=0$, cf.\ \eqref{muu}. Then for $\nu$-almost any $z \in \widehat{\Sigma}$ and any fixed $\omega >0$, we have
\begin{equation}\label{f1}
u_{\mu}/u_{\nu}(x) \ra d \mu/d \nu(z)\mm{ as } x \ra z, \mm{ with } x \in \P(z,\omega).
\end{equation}
where $\P(z,\omega):= \{x \in H \setminus \Sigma \,|\, \delta_{\bp} (x) >\omega \cdot d_H(x,z)\}$ non-tangential \si-pencils and $d \mu/d \nu$ is the Radon-Nikodym derivative of $\mu$ with respect to $\nu$.
\end{proposition}
For $H \in {\cal{H}}^{\R}_n$ we also have a point at infinity $\infty_H \in \widehat{\Sigma}$, in this case we choose some basepoint $p_0$ and the version of (\ref{f1}) for $\infty_H$  reads
\begin{equation}\label{f2}
u_{\mu}/u_{\nu}(x) \ra d \mu/d \nu(\infty)\mm{ as } x \ra \infty, \mm{ with } x \in \P(p_0,\omega).
\end{equation}

\noindent\textbf{Proof} \,   \cite{L2},Theorem 4 (and Lemma 4.8 for the result in $\infty_H \in \widehat{\Sigma}$). \qed

Note that in the case of a minimal solution the associated measure $\nu$ is a Dirac measure concentrated in one point $p$ and the latter result says that, for any fixed $\omega >0$ we always have $u_{\mu}/u_{\nu}(x) \ra d \mu/d \nu(z)\mm{ as } x \ra p, \mm{ with } x \in \P(p,\omega)$.

\setcounter{section}{3}
\renewcommand{\thesubsection}{\thesection}
\subsection{Minimal Hypersurfaces in $E^\pm_C$} \label{i3}

In this chapter we show the asserted existence of the hypersurface $H_C \subset E_C$  and prove its regularity under the assumption, we settle in the next chapter, that the collection $\tau \cdot H_C \subset E_C$, for $\tau \in (0,1)$ of scaled copies of  $H_C \subset E_C$ consists of disjoint leaves.

\subsubsection{Existence Results} \label{r1}

We recall from \cite{BG} that any  oriented minimal boundary $H \subset \R^{n+1}$ is connected and its complement has exactly two path components $E_H^+ \cup E_H^-$. Henceforth $E_C$ will always denote one of the two path components $E_C^+ \cup E_C^-=\R^{n+1} \setminus C$. More generally, we denote by $E_H$ one of the two components $E_H^+ \cup E_H^-$ in the complement of an oriented minimal boundary $H \subset \R^{n+1}$.

\begin{lemma}\label{ti}   \, Let $H^n \subset \R^{n+1}$ be an oriented minimal boundary. Then there is a (generally non-unique) area minimizing cone $C$ approximating $H$ near infinity called a \textbf{tangent cone at infinity}. Then one can see that for sufficiently small $\ve >0$, $H$ can locally be written as a smooth section of the normal bundle of $C$ over $\bp_C^{-1}(\ve) \subset C$.
\end{lemma}

 \textbf{Proof} \, The Allard type graphical representation and the proof are the same as for tangent cones around singular points \cite{Gi}, Ch.9, with the one difference that we scale $H$ by ever smaller constants $\tau \ra 0$.\qed

\begin{lemma}\label{fo}\, Let $C^n \subset \R^{n+1}$ be an area minimizing cone. Then there is an oriented minimal boundary $H_C \subset E_C$ and  $C$  is its unique tangent cone at infinity.
\end{lemma}

The construction of $H_C$ essentially carries over from \cite{HS}, p.113-114 where $C$ was assumed to be singular only in $0$. The only difference to \cite{HS}  is that we additionally use the strict maximum principle \cite{Si}. We include the main steps of  the argument for the reader's convenience.\\

 \textbf{Proof} \, Choose area minimizing Plateau solutions $P_\ve \subset \overline{E_C \cap B_1(0)}$, bounding an open set in $B_1(0)$ relative $\p B_1(0)$ with a smooth prescribed border $\p P_\ve=\p B_1(0) \cap P_\ve$ and $\p P_\ve \cap (\p B_1(0) \cap C)\v$  so that the Hausdorff-distance of $\p P_\ve$ to $\p B_1(0) \cap C$ is at most $\ve>0$. The condition $\p P_\ve \cap (\p B_1(0) \cap C) \v$ implies  $P_\ve \cap (\overline{B_1(0)} \cap C) \v$ from the use of the strict maximum principle needed since $C \setminus \{0\}$ may be singular.\\

 In turn, one finds that for $\ve \ra 0$ the $P_\ve$ converge to some area minimizing Plateau solution $P_0$ in $\overline{B_1(0)}$ with  $\p P_0   = \p B_1(0) \cap C$ since both, $C$ and $C \setminus B_1(0) \cup P_0$ are area minimizing and singularities must have codimension $\ge 7$ one infers that  $P_0= \overline{B_1(0)} \cap C$. Thus $dist(0,P_\ve) \ra 0$ for $\ve \ra 0$. Then $dist(0,P_\ve)^{-1} \cdot P_\ve$ subconverges to an oriented minimal boundary $H \subset E_C$ with $dist(0,H) = 1$. We observe that $\tau_i  \cdot H$, for any sequence $\tau_i \ra 0$, converges to $C$ since the limit hypersurface must be scaling invariant. This renders
 $C$ as the unique tangent cone of $H$ at infinity.  \qed

\subsubsection{Jacobi fields and Variations of Area Minimizers} \label{je}

We discuss two methods to associate a Jacobi field to the area minimizer $H_C$. We start with the classical smooth variations  $H_t$, $t \in [0,1]$,  of a given area minimizer $H_0=H$ through a family of others gives rise to Jacobi fields along $H$ from the derivative in $t$-direction. In our case, one is tempted to use $H_0=C$
and $H_\tau:= \tau \cdot  H_C$, $\tau  \in (0,1]$, to define such a variation. We observe that since the $\tau \cdot  H_C$ are all on one side of $C$ the normal component of the induced Jacobi field can be written as $f \cdot \nu$, for some $f \ge 0$, where $\nu$ is unit normal vector field of $C$ pointing into $E_C$.  This  clearly only applies to the smooth portions $C \setminus \Sigma_C$ of $C$  where we can locally write $\tau \cdot  H_C$ as a smooth section of $\nu$, when $\tau >0$ becomes small enough.\\

A byproduct if the fact that  $C$ is the tangent cone of $H_C$ at infinity is that $H_C$ typically decays so fast towards $C$ that the Jacobi field defined from this parametrization of the variational family simply vanishes. \\

Thus we seek for a suitable $\tau$-reparametrization of $H_\tau$, $\tau \in (0,1]$, to get a non-vanishing Jacobi field. To this end the reparametrized family must also be differentiable for $\tau=0$, in all points of  $C \setminus \Sigma_C$. However, for an arbitrary resulting Jacobi field this can only be granted in one single point. We accomplish the extension to all points from an a priori understand of the possibly resulting Jacobi fields.We show inductively that the only options are Jacobi fields $f \cdot \nu$ where either $f$ is a convex combination of the two minimal solutions associated to the Dirac measure in $0$ or in $\infty$\\

To reach this a priori insight we use another way to extract Jacobi fields from the area minimizer $H_C$ and now we actually employ the fact that  $H_C$ decays towards $C$. For small $\ve >0$ we find that the $\tau \cdot H_C$, for $\tau \in (0,1]$, can be written as section $u_\tau \cdot \nu$ of the normal bundle $\nu$ over $\bp^{-1}(\ve)\subset C$. The minimal surface equation for such graphs has the form
\begin{equation}\label{nlm}
\Delta \, u_\tau + |A|^2 \cdot u_\tau = div(a \cdot \nabla u_\tau) + b \cdot \nabla u_\tau + c \cdot u_\tau
\end{equation}
where, however, the coefficients also depend on $u_\tau$, e.g.\cite{Si},p.333,eq.(7). The point is that these coefficients $C^k$-compactly converging to zero, for any $k \in \Z^{\ge 0}$, for the corresponding convergence of $ u_\tau$ to zero. Thus we can fix some basepoint $p \in \bp^{-1}(\ve)$ and consider the functions $u_\tau/u_\tau(p)$. For $\tau \ra 0$, the Harnack inequality ensures that we get subconverging sequences with limits
$v>0$ solving $\Delta \, v + |A|^2 \cdot v = 0$  on $\bp^{-1}(\ve)$  and due to scaling process this can be iterated to give a Jacobi field on $C \setminus \Sigma_C$.\\

\subsubsection{Foliations and Regularization} \label{r2}

Our main goal is to show  for any area minimizing cone $C^n \subset \R^{n+1}$ any oriented minimal boundary $H_C \subset E_C$ must be regular.  The idea of \cite{HS} in the case of $C$  singular only in $0$  was to represent $H$ as a polar graph and this suffices to prove its regularity. This strategy and the tools used in \cite{HS}  do not apply in the general case. In our approach the main challenge is to show that $H \cap \ve \cdot H \v$ for some suitably small $\ve >0$. This is what will occupy us in the next chapter. From this we  will be derive that $H \cap s\cdot H \v$ for any $s \in (0,1)$ and this is enough to show that $H$ is regular since we prove the following regularity result:

 \begin{proposition}\label{ceg}    \, Let $H \subset \R^{n+1}$ be an oriented minimal boundary in $\R^{n+1}$ and assume that $s\cdot H \cap t \cdot H \v$, for any two $s \neq t \in (0,1)$, then $H$ is regular.
\end{proposition}

 \textbf{Proof} \,  Assume that $H$ is singular in some $x \in H$. Then we choose some singular tangent cone $C$ in the singular point $x$ and think of it as a hypersurface in $\R^{n+1}$ with basepoint $=$ tip in $x$. Then $s \cdot C$ becomes a tangent cone
 of $ s\cdot H$ in the singular point $s \cdot x$ for any $s \in (0,1)$. Since $H$ is an oriented boundary $H = \p Y$ for some open set $Y\subset \R^{n+1}$, which we can choose so that  $s\cdot H \subset Y$ for $s \in (0,1)$ and, hence, $s\cdot H \subset  t \cdot Y$, for
  $s<t$, $s,t \in (0,1)$. We choose $s=1/2$ and $t_i=1/2+1/2^i$, for $i \in \Z^{>0}$, and apply two linear transformations to these hypersurfaces and the tangent cones. We scale the configuration by  $a_i:=d(s \cdot x, t_i \cdot x)^{-1}$ to gauge the distance between the two singular points to 1. Secondly, we translate the configuration by $-a_i \cdot s \cdot x$. The two transformed singular points in the two hypersurface  $H_{0,i}:= a_i \cdot s \cdot H -a_i \cdot s \cdot x$ and  $H_{i,i}:= a_i \cdot t_i \cdot H -a_i \cdot s \cdot x$
are $0$ and $v:=x/|x|$.\\

 Since $a_i \ra \infty$, for $i \ra \infty$  we infer, using  $s\cdot H \subset  t \cdot Y$, for $s<t$, that $C+v \subset \overline{E_C}$, where $E$ is one of the two path components of $\R^{n+1} \setminus C$. The strict maximum principle even shows that $C+v\subset E_C$. Since $\tau \cdot C=C$ and $\tau \cdot C_v= C_{\tau \cdot v}$ we see $C \cap C_{\tau \cdot v} \v$, for any $\tau \in (0,1)$. The  component $E_C$ also contains the  $C_{\tau \cdot v}$, $\tau \in (0,1)$. The family $C_{\tau \cdot v}$, $\tau \in [0,1]$ is a variation of $C$ with all hypersurfaces situated on one side of $C$. It defines a nowhere vanishing Jacobi field $\alpha_C$ on $C \setminus \Sigma_C$. When we consider any of the regular rays $\R^{>0} \cdot \nu \in C \setminus \Sigma_C$, for some $\nu \in  \p B_1(0) \cap C$, we observe that  $\alpha_C$ is constant along $\R^{>0} \cdot \nu$, that is, its is parallel with respect to the covariant derivative inherited from $\R^{n+1}$, and the same holds for its normal component relative $C$. However, since the normal component of $\alpha_C$ is a positive solution $J_C (\alpha_C)=0$, we know from \ref{mc} that there is no solution constant along $\R^{>0} \cdot \nu$ unless $C$ was a hyperplane. \qed

\setcounter{section}{4}
\renewcommand{\thesubsection}{\thesection}
\subsection{Separation of Hypersurface Leaves} \label{4}

Now we reach the central interplay of minimal growth properties of Jacobi fields and associated area minimizing variations of the cone.
The arguments are by induction the maximal number ${\cal{N}}(C)$ of iterated tangent cone blow-ups, outside the respective origin of the tangent cone, until we reach exclusively tangent cones which can be written as products of cones with isolated singularities with some $\R^k$. We observe that ${\cal{N}}(C)$ can be any integer between $n-7$ and $0$.

\subsubsection{Isolated Singularities ${\cal{N}}(C)=0$} \label{s1}
We start the induction scheme with the case of a cone $C$ singular only in $0$, as treated in \cite{HS}. However, we give a new proof since we want to derive a stronger result.
\begin{itemize}
  \item For $\lambda^{\bp}_{J,C} = 0$ the only positive solution of $J_C \, w=0$ is $\Psi(\omega,r) = \psi(\omega) \cdot r^{- \frac{n-2}{2}}$.
  \item For $\lambda^{\bp}_{J,C} > 0$ all positive solutions are
of the form \[a \cdot \Psi_+(\omega,r)  + b \cdot \Psi_-(\omega,r)  = \psi(\omega) \cdot (a \cdot r^{\alpha_+} + b \cdot \cdot r^{\alpha_-}),\] for some $a,b \ge 0, a+b >0$ and $\alpha_- < \alpha_+<0$.
\end{itemize}
Now let $H \subset E_C$ be an oriented minimal boundary with unique tangent cone $C$ at infinity. From the discussion in the last section of Ch.\ref{je} we infer there is a radius $\rho>0$ so that for any ray $\R^{>0} \cdot \eta \in C$, for some $\eta \in  \p B_1(0) \cap C$, $H|_{\R^{>\rho} \cdot \eta}$ is a \textbf{strictly decreasing} function. This means that $(H \cap \tau \cdot H) \cap \R^{n+1} \setminus B_\rho(0) \v$, $\tau \in (0,1)$ and from the strict maximum principle we therefore have $H \cap \tau \cdot H \v$, $\tau \in (0,1)$. Thus \ref{ceg} implies that $H$ is a \textbf{regular} hypersurface.\\

From \ref{mc} we know that for  $\lambda^{\bp}_{J,H} > 0$ its Martin boundary is just the point $\{\infty\}$. For $\lambda^{\bp}_{J,H} > 0$ we only have one positive solution, the ground state.
We know from \ref{ver} that in both cases the convergence $\tau \cdot H \ra C$, for $\tau \ra 0$ only induces solutions on $C$ with \textbf{minimal growth} of towards $\Sigma_C$.\\

To prove the \textbf{uniqueness} assertion we start with some oriented minimal boundary $T \subset E_C$ and we observe that its tangent cone at infinity $C_T$ satisfies  $C_T \subset \overline{E_C}$ and $C_T \cap C \supset \{0\} \n$. Thus the strict maximum principle shows that $C_T=C$. Then we infer as in the case of $H$, that $T$ is regular and we get the same Jacobi field growth estimates as for $H$. We may assume there are constants $c_u \ge c_l>0$ so that $E_{c_l \cdot H} \subset E_T \subset  E_{c_u \cdot H}$. (Otherwise we could build an area minimizing hypersurface $Z$ with $Z \cap B_r(0) = b \cdot T \cap B_a(0)$ and $Z \setminus B_b(0) = b \cdot H \setminus B_R(0)$, for some $a, b>0$ and $R>r>0$ with singularities of codimension $\le 1$.) For sufficiently large $i$ we see that on $B_{2^{i+1}}(0) \setminus B_{2^i}(0) \cap C$ the constant $c_u/c_l$ can be chosen arbitrarily near to $1$ and from the strict maximum principle also on $B_{2^i}(0) \cap C$. For $i \ra \infty$, we see that $H= k \cdot T$, for some $k >0$. This concludes the proof of our main theorem for ${\cal{N}}(C)=0$. \qed

\subsubsection{Inductive Arguments} \label{ia}

In this section we assume we established the main theorem for all area minimizing cones with ${\cal{N}}(C) \le m$ and use this to derive the result for cones with ${\cal{N}}\le m+1$.

\begin{proposition}\label{ti}   \, Let $C$ be a cone with ${\cal{N}}(C) \le m$ with minimal smoothing $H_{C}$. Then the smooth hypersurface $H_{C}\times \R \subset E_{{C} \times \R}$ is again the unique such area minimizer in $E_{{C} \times \R}$ and it induces the unique Jacobi field on ${C} \times \R$ with minimal growth towards $\Sigma_C \times \R$.
\end{proposition}

 \textbf{Proof} \, We know that $H \times \R$ corresponds to the Jacobi field $J$ represented by the unique Martin boundary point at infinity. It has minimal growth towards $\Sigma$. Given another area minimizing hypersurface $F \subset  E_{C \times \R}$. There are two possible cases for the Martin measure $\mu_u$ of the associated Jacobi field  $u=\int_{\Sigma_C \cup \{\infty\}}k(x;y) \, d \mu_u$. Either  $\mu_u$ has vanishing or positive Radon-Nikodym derivative relative the Dirac measure $\delta_\infty$ associated to $J$. In both cases we can use the Fatou theorem \ref{ft} which shows that in the first case $u$ has a stronger decay towards infinity than $J$. But then we can scale $F$ so that it intersects $\R \times H$  in non-empty compact set, using the result  \ref{ti} that follows. As in Ch.\ref{s1} we get a contradiction to regularity of area minimizers.\\

Thus, $u/J$ converges non-tangentially to a positive constant, the Radon-Nikodym derivative and we get from \ref{ti} that for some $k \ge 1$:  $F$ lies between   $k \cdot (H \times \R)$ and  $C \times \R$
and $H \times \R$ lies between $k \cdot F$  and  $C \times \R$. As in Ch.\ref{s1} we conclude that $H \times \R = c \cdot F$, for some $c>0$.\qed

\begin{lemma}\label{ti}   \, Let $C$ be an area minimizing cone and assume there are two area minimizing hypersurfaces   $H_{C}, H^*_{C} \subset E_{C} $, so that both of them admit local smooth graphical representations $f_H \cdot \nu_C$ and  $f_{H^*} \cdot \nu_C$ as sections of the normal bundle over $\bp^{-1}(1)$. If there a constant $c(H,H^*) \ge 1$, so that if $f_H > c \cdot  f_{H^*}$ over $\bp^{-1}(1)$, then, for sufficiently small $d>0$, we have   $d \cdot H^* \cap H \v$.
\end{lemma}

 \textbf{Proof} \, Otherwise, we choose a sequence of $d_i \ra 0$ and consider the closed set $\overline{d_i \cdot H^* \cap E_H}$. The Harnack inequality for $J$  shows that we may assume that $f_H > c \cdot  f_{H^*}$ holds beyond $\bp^{-1}(1)$ on any  subcone $C(A)= A \times \R^{> 0} \subset \p B_1(0) \cap C \times \R^{> 0} =C \setminus \{0\}$, for some relative compact open set $A \in \p B_1(0) \cap C$, outside suitably large balls. It cannot contain any compact components as is seen from the same regularity argument as in Ch.\ref{s1}. This argument further urges a superpolynomial growth of $Area(B_r(0) \cap \overline{d_i \cdot H^* \cap E_H})$, for $r \ra \infty$. Since $H$ and $H^*$ are oriented boundaries, we find can construct a subconverging sequence, for  $d_i \ra 0$, and get a contradiction to the non-extinction lemma A.7 of \cite{L1}.\qed

 Now we study the general case cones and $H_C \subset E_C$ with ${\cal{N}}= m+1$. We show, in two steps, that the associated Jacobi field has minimal growth towards $\Sigma_C$. We first use an induction argument to show this for $\Sigma_C \setminus \{0\}$. From this we infer the smoothness of $H_C \subset E_C$ and, as in the case of isolated singularises we infer the minimal growth towards $\{0\}$.   Finally, to conclude the induction step, we prove the uniqueness assertion.

 \begin{proposition}\label{jg0}   \, For area minimizing cone $C$ with ${\cal{N}}= m+1$ we consider an area minimizing hypersurface $H_C \subset E_C$ and any associated Jacobi field $J$.
 Then $J$ has minimal growth towards $\Sigma_C \setminus \{0\}$
\end{proposition}

 \textbf{Proof} \, As in the case of isolated singularities any tangent cone  of $T$ at infinity $C_T$ satisfies  $C_T \subset \overline{E_C}$ and $C_T \cap C \supset \{0\} \neq$ and, thus, $C_T=C$ from the strict maximum principle. Now we observe that $T$ induces further area minimizing hypersurfaces $H_{C_p} \subset E_{C_p}$ for each of its tangent cones $C_p$, in any $p \in \Sigma \setminus \{0\}$. All of these tangent cones are of the form $C^* \times \R$, for some area minimizing cone $C^*$. Thus, form the induction hypothesis applied to \ref{ti} we infer that $H_{C_p}$ is smooth and $H_{C_p}=H_{C^*}\times \R \subset E_{{C^*} \times \R}$ up to multiples. We also know from \ref{ti} that $H_{C_p}$ induces a Jacobi field of minimal growth towards $\Sigma_{C_p} \times \R$. In turn, \ref{mc} and \ref{tri} show that this function is again uniquely determined. Thus it is the solution induced from any Jacobi field $J$ associated to  $H_C \subset E_C$. But this means, from \ref{ver}, that the  $J$  has minimal growth towards $\Sigma_C \setminus \{0\}$. \qed

As a consequence of \ref{mc} and \ref{tri}, now applied to $C$, there are only tow possible cases for $J$ left:
\begin{itemize}
  \item For $\lambda^{\bp}_{J,C} = 0$ the only positive solution of $J_C \, w=0$ is $\Psi(\omega,r) = \psi(\omega) \cdot r^{- \frac{n-2}{2}}$.
  \item For $\lambda^{\bp}_{J,C} > 0$ all positive solutions are
of the form \[a \cdot \Psi_+(\omega,r)  + b \cdot \Psi_-(\omega,r)  = \psi(\omega) \cdot (a \cdot r^{\alpha_+} + b \cdot \cdot r^{\alpha_-}),\] for some $a,b \ge 0, a+b >0$ and $\alpha_- < \alpha_+<0$.
\end{itemize}

Now let $H \subset E_C$ be an oriented minimal boundary with unique tangent cone $C$ at infinity. We infer that on any  subcone $C(A)= A \times \R^{> 0} \subset \p B_1(0) \cap C \times \R^{> 0} =C \setminus \{0\}$, for some relative compact open set $A \in \p B_1(0) \cap C$, and outside some suitably large ball for any ray $\R^{>0} \cdot \nu \in C$, for some $\eta \in  A$, $H|_{\R^{>\rho} \cdot \eta}$ is a strictly decreasing function. This means that $(H \cap \tau \cdot H) \cap \R^{n+1} \setminus B_\rho(0) \v$, $\tau \in (0,1)$. From \ref{ti} we can again conclude that $H \cap \tau \cdot H \v$, $\tau \in (0,1)$. Thus \ref{ceg} implies that $H$ is a regular hypersurface. From \ref{mc} we know that for  $\lambda^{\bp}_{J,H} > 0$ its Martin boundary is just the point $\{\infty\}$. For $\lambda^{\bp}_{J,H} > 0$ we only have one positive solution, the ground state.
We know from \ref{ver} that in both cases the convergence $\tau \cdot H \ra C$, for $\tau \ra 0$ exclusively induces solutions on $C$ with  minimal growth of towards $\Sigma_C$.\\

To prove the  uniqueness assertion we start with some oriented minimal boundary $T \subset E_C$. We have already seen,  in the proof of \ref{jg0} above, that $C_T=C$. We infer as in the case of $H$, that $T$ is regular and we get the same Jacobi field growth estimates as for $H$. From \ref{ti} we may assume there are constants $c_u \ge c_l>0$ so that $E_{c_l \cdot H} \subset E_T \subset  E_{c_u \cdot H}$. (Otherwise we could build an area minimizing hypersurface $Z$ with $Z \cap B_r(0) = b \cdot T \cap B_a(0)$ and $Z \setminus B_b(0) = b \cdot H \setminus B_R(0)$, for some $a, b>0$ and $R>r>0$ with singularities of codimension $\le 1$.) For sufficiently large $i$ we see that on $B_{2^{i+1}}(0) \setminus B_{2^i}(0) \cap C$ the constant $c_u/c_l$ can be chosen arbitrarily near to $1$ and from the strict maximum principle also on $B_{2^i}(0) \cap C$. For $i \ra \infty$, we see that $H= k \cdot T$. \qed

\end{document}